\documentclass[a4paper,12pt]{amsart}
\usepackage{graphicx}
\usepackage{amscd}
\usepackage{amssymb}
\newtheorem{thm}{Theorem}[section]

\newtheorem{lem}[thm]{Lemma}
\newtheorem{prop}[thm]{Proposition}
\theoremstyle{defn}
\newtheorem{defn}[thm]{Definition}
\theoremstyle{rem}

\numberwithin{equation}{section}

\newcommand{\eps}{\varepsilon}

\newcommand{\G}{\mathcal{G}}

\newcommand{\T}{\mathcal{T}(\mathcal G)}
\newcommand{\U}{\mathcal{U}}
\newcommand{\E}{\mathcal{E}nd(\mathcal U)}

\newcommand{\Euv}{\mathcal E_{u,v}^\pi}
\newcommand{\R}{\mathcal{R}ep(\mathcal G)}
\newcommand{\F}{\mathfrak {F}_{u,v}}
\newcommand{\HH}{\mathcal H_\pi}
\newcommand{\Hu}{\mathcal H_u^\pi}
\newcommand{\Hv}{\mathcal H_v^\pi}

\newcommand{\repn}{representation\,}
\newcommand{\repns}{representations\,}
\newcommand{\nt}{natural transformation\,}

\newcommand{\oo}{\otimes}
\begin{document}

\title[duality for groupoids]{Tannaka-Krein duality for compact groupoids III,
duality theory}%
\author{Massoud Amini}%
\address{
Department of Mathematics, Tarbiat Modarres University, P.O.Box
14115-175, Tehran , Iran , mamini@modares.ac.ir \newline
Department of Mathematics and Statistics, University of
Saskatchewan, 106 Wiggins Road, Saskatoon, S7N 5E6 ,
mamini@math.usask.ca}

\thanks{The author was visiting the University of Saskatchewan during the
preparation of this work. He would like to thank University of
Saskatchewan and Professor Mahmood Khoshkam for their hospitality
and support}%

\begin{abstract}
In a series of papers, we have shown that from the \repn theory of
a compact groupoid one can reconstruct the groupoid using the
procedure similar to the Tannaka-Krein duality for compact groups.
In this part we introduce the Tannaka groupoid of a compact
groupoid and show how to recover the original groupoid from its
Tannaka groupoid.
\end{abstract}
\maketitle
\section{introduction}
This is the last in a series of papers in which we generalized the
Tannaka-Krein duality to compact groupoids. In [A1] we studied the
\repn theory of compact groupoids. In particular, we showed that
irreducible \repns have finite dimensional fibres. We also proved
the Schur's lemma, Gelfand-Raikov theorem and Peter-Weyl theorem
for compact groupoids. In [A2] we studied the Fourier and
Fourier-Plancherel transforms and their inverse transforms on
compact groupoids. In this part we show how to recover a compact
groupoid from its \repn theory. This is done along the lines of
the Tannaka duality for compact groups. We refer the interested
reader to [JS] for a clear exposition of this theory. All over
this paper we assume that $\G$ is compact and the Haar system on
$\G$ is normalized. We put $X=\G^{(0)}$.

\section{Tannaka groupoid}

There is a forgetful functor $\U:\R\to \mathcal{H}il_X$ to the
category of Hilbert bundles over $X$ and  operator bundles. A {\it
natural transformation} $a:\U\to\U$ is a family of bundle maps
$a_\pi:\HH\to\HH$ indexed by $\R$ such that for each
$\pi_1,\pi_2\in\R$ and $h\in Mor(\pi_1,\pi_2)$ the following
diagram commutes

\begin{equation*}
\begin{CD}
\mathcal H_{\pi_1}@>{a_{\pi_1}}>>\mathcal H_{\pi_1}\\
@VhVV    @VVhV\\
\mathcal H_{\pi_2}@>>{a_{\pi_2}}>\mathcal H_{\pi_2}
\end{CD}
\end{equation*}

One should understand this as each $a_\pi$ being a bundle
$a_\pi=\{a_{u,v}^\pi\}$ of bounded linear operators
$a_{u,v}^\pi\in\mathcal B(H_u^{\pi},\mathcal H_v^{\pi})$ (possibly
zero) indexed by $X\times X$ such that for each $u,v\in X$ the
following diagrams commute

\begin{equation*}
\begin{CD}
\mathcal H_{u}^{\pi_1}@>{a_{u,v}^{\pi_1}}>>\mathcal H_{v}^{\pi_1}\\
@V{h_{u}}VV    @VV{h_{v}}V\\
\mathcal H_{u}^{\pi_2}@>>{a_{u,v}^{\pi_2}}>\mathcal H_{v}^{\pi_2}
\end{CD}
\end{equation*}

Given $x\in\G$ there is a natural transformation $\mathcal
T_x:\U\to\U$ defined by
\begin{equation}
(\mathcal T_x)_{u,v}^\pi=
\begin{cases}
\,\pi(x)& \quad \text{if} \,\, u=s(x),v=r(x),\\
\,0&\quad \text{otherwise}
\end{cases}
\end{equation}

Another interesting example of a natural  transformation is the
Fourier transform [A2]. Recall that we looked at $L^1(\G)$ as a
bundle of Banach algebras over $\G^{(0)}\times\G^{(0)}$ whose
fiber at $(u,v)$ is $L^1(\G_u^v,\lambda_u^v)$, and then each $f\in
L^1(\G)$ had its Fourier transform $\mathfrak F(f)$ in $C_0(\hat
\G,\mathcal B(\mathcal H))$, where $\mathcal B(\mathcal H)$ is a
bundle of bundles of $C^{*}$-algebras over $\hat \G$ whose fiber
at $\pi$ is the bundle $\mathcal B(\HH)$ over
$\G^{(0)}\times\G^{(0)}$ whose fiber at $(u,v)$ is $\mathcal
B(\Hv,\Hu)$, the space $C_0(\hat G,\mathcal B(\mathcal H))$ is the
set of all continuous sections vanishing at infinity, and
$\mathfrak {F}(f)(\pi)_{(u,v)}=\F(f_{(u,v)})(\pi)$. Now we need a
flip in the order of $u,v$ when we consider $\mathfrak F(f)$ as a
natural transformation, namely we put $\mathfrak F(f)_{u,v}^\pi=
\mathfrak {F}(f)(\pi)_{(v,u)}$. This way we get $\mathfrak
F(f)_{u,v}^\pi\in\mathcal B(\Hu,\Hv)$. To see that this is indeed
a natural transformation, note that for each $u,v\in X$,
$x\in\G_v^u$, $\pi_1,\pi_2\in\R$, and $h\in Mor(\pi_1,\pi_2)$ we
have

\begin{equation*}
\begin{CD}
\mathcal H_{u}^{\pi_1}@>{\pi_1(x^{-1})}>>\mathcal H_{v}^{\pi_1}\\
@V{h_{u}}VV    @VV{h_v}V\\
\mathcal H_{v}^{\pi_2}@>>{\pi_2(x^{-1})}>\mathcal H_{v}^{\pi_2}
\end{CD}
\end{equation*}

Multiplying both sides with $f_{(v,u)}(x)$ and integrating against
$d\lambda_v^u(x)$ we get

\begin{equation*}
\begin{CD}
\mathcal H_{u}^{\pi_1}@>{\mathfrak F(f)_{u,v}^{\pi_1}}>>\mathcal H_{v}^{\pi_1}\\
@V{h_{u}}VV    @VV{h_v}V\\
\mathcal H_{v}^{\pi_2}@>>{\mathfrak F(f)_{u,v}^{\pi_2}}>\mathcal
H_{v}^{\pi_2}
\end{CD}
\end{equation*}

which means $\mathfrak F(f):\mathcal U\to\mathcal U$ is a natural
transformation. Let $\mathcal{E}nd(\U)$ be the set of all natural
transformations $:\U\to\U$ with the coarsest topology making all
maps $a\mapsto a_{u,v}^\pi$ continuous. Also we define an
involution on $\E$ by
$$\bar a_{u,v}^\pi(\xi)=\overline{a_{u,v}^{\bar\pi}(\bar\xi)}\quad
(u,v\in X,\pi\in\R,\xi\in\overline{\mathcal H_u^\pi}).$$ The
following is trivial.

\begin{lem}
$\E$ is a topological vector space with continuous involution.\qed
\end{lem}

\begin{prop}
The map $$q:\E\to \prod_{\rho\in\hat\G} \mathcal{E}nd(\mathcal
H_\rho)$$
$$a\mapsto (a_\rho)_{\rho\in\hat\G}$$
is an isomorphism of topological vector spaces.
\end{prop}
{\bf Proof} The following commutative diagrams (with  vertical
maps being canonical imbedding ) illustrates that
$a_{{\pi_1}\oplus{\pi_2}}=a_{\pi_1}\oplus a_{\pi_2}$, for each
$\pi_1,\pi_2\in\R$ and $a\in\E$.

\begin{equation*}
\begin{CD}
\mathcal H_{\pi_1}@>{a_{\pi_1}}>>\mathcal H_{\pi_1}\\
@V{\iota_1}VV    @VV{\iota_1}V\\
\mathcal H_{\pi_1}\bigoplus\mathcal H_{\pi_2}@>>
{a_{\pi_1\oplus\pi_2}}>\mathcal H_{\pi_1}\bigoplus\mathcal H_{\pi_2}\\
@A{\iota_2}AA    @AA{\iota_2}A\\
\mathcal H_{\pi_2}@>>{a_{\pi_2}}>\mathcal H_{\pi_2}
\end{CD}
\end{equation*}
This  plus the fact that each \repn of $\G$ is the direct sum of
its irreducible sub representations [A1, theorem 2.16] shows that
$q$ is one-one. To show that it is onto, let $b=(b_\rho)$ with
$b_\rho\in\mathcal{E}nd(\mathcal H_\rho)$ be given. Let $\pi\in\R$
and $\mathcal H_\pi=\bigoplus_{\rho\in\hat\G} \mathcal
H_{\pi_\rho}$ be the unique decomposition into isotropical
components. For $\rho\in\hat\G$, the canonical map
$$\psi_\rho:\mathcal H_\rho\bigotimes
Hom_\G(\mathcal H_\rho,\mathcal H_\rho)\to\mathcal H_\rho$$
$$\xi\otimes\varphi\mapsto \varphi(\xi)$$
is  an isomorphism of $\G$-modules. Put
$a_{\pi_\rho}=\psi_\rho\circ(b_\rho\otimes id)
\circ\psi_\rho^{-1}$ , $a_\pi=\bigoplus_{\rho\in\hat\G}
a_{\pi_\rho}$, and $a=(a_\pi)_{\pi\in\R}$. It is easy to see that
$a:\U\to\U$ is a \nt and $a_\rho=b_\rho$, for each
$\rho\in\hat\G$. Hence $q$ is onto. The way we defined the
topology of $\E$ makes $q^{-1}$ continuous. The fact that $q$ is
continuous is trivial.\qed

An  alternative version of the above proposition would be to
interpret $q$ as a bundle isomorphism between bundles of bundles
of $C^*$-algebras, that is $\E$ is a bundle over $\hat \G$ whose
fiber at $\rho\in\hat \G$ is a bundle of $C^*$-algebras over
$X\times X$ whose fiber at $(u,v)$ is the $C^*$-algebra $\mathcal
B(\Hu,\Hv)$. This has the advantage of a better interpretation of
the global Fourier transform. Indeed, in the light of [A2,
Corollary 2.4], the above discussion could be rephrased as

\begin{prop}
The global Fourier transform
$$\mathfrak F:L^1(\G)\to\E$$
is a bundle homomorphism.\qed
\end{prop}

Although $\E$ doesn't seem to have a non  trivial  everywhere
defined product, but one can define a "center" for it!

\begin{defn}
The center $\mathcal Z(\E)$ of $\E$ consists  of those $a\in\E$
which commute with each $b\in\E$ in the following sense
$$b_{v,u}^\pi\circ a_{u,v}^\pi=a_{v,u}^\pi\circ b_{u,v}^\pi\quad(u,v\in X,\pi\in\R).$$
\end{defn}

\begin{prop}
$\mathcal Z(\E)$ is a closed subspace of $\E$  and the restriction
of $q$ gives an isomorphism
$$\mathcal Z(\E)\simeq \mathbb C^{\hat \G}=\prod_{\rho\in\hat\G} \mathbb C.id_\rho.$$
\end{prop}
{\bf Proof} The first statement is trivial. The second follows by
a  diagonalization argument.\qed

\vspace{.3 cm} Using the notion of center, some of the results of
the previous section on the Fourier transform of central functions
could be rephrased in the terms of $\E$. Here is one example.

\begin{lem}
Let $a\in\E$,  then $a\in\mathcal Z(\E)$ if and only if there is
$\pi\in\R$ such that $a_\pi:\HH\to\HH$ commutes with the action of
$\G$, i.e.
$$a_{v,u}^\pi\pi(x)=\pi(x^{-1})a_{u,v}^\pi\quad(u,v\in X, x\in\G_u^v).$$
\end{lem}
{\bf Proof} If $a\in\mathcal Z(\E)$, then for  each $x\in\G$, $a$
commutes with $\mathcal T^x$, so we have the above equality.
Conversely, if this holds, then for each $\pi\in\R$, $a_\pi\in
Mor(\pi,\pi)$, so by the definition of the natural transformation,
$a$ commutes with each $b\in\E$.\qed

Now Lemma 4.3 of [A2] could be rephrased as

\begin{prop}
If $f\in\mathfrak CC(\G)$ then $\mathfrak {DF}(f)\in\mathcal
Z(\E)$.\qed
\end{prop}

\begin{defn}
An element $a\in\E$ is called monoidal (tensor preserving) if for
each $\pi_1,\pi_2\in\R$,  $a_{\pi_1\oo\pi_2}=a_{\pi_1}\oo
a_{\pi_2}$ and $a_{tr}$ is trivial, i.e. for each $u,v\in X$ the
following diagram commutes

\begin{equation*}
\begin{CD}
\mathcal H_{u}^{\pi_1}\bigotimes\mathcal H_{u}^{\pi_2}
@>{a_{u,v}^{\pi_1}\oo a_{u,v}^{\pi_2}}>>\mathcal
H_{v}^{\pi_1}\bigotimes
\mathcal H_{v}^{\pi_2}\\
@|    @|\\
\mathcal
H_{u}^{\pi_1\oo\pi_2}@>>{a_{u,v}^{\pi_1\oo\pi_2}}>\mathcal
H_{v}^{\pi_1\oo\pi_2}
\end{CD}
\end{equation*}

and $a_{u,v}^{tr}=id$ , where $tr$ is the trivial \repn of $\G$ on
$\mathbb C$. $a\in\E$ is called Hermitian if $\bar a=a$.
\end{defn}

\begin{defn}
For each $u,v\in X$ and $a\in\E$, consider the continuous section
$a_{u,v}$ defined on $\R$ by $a_{u,v}(\pi)=a_{u,v}^\pi$. The set
$\T$ of all sections $a_{u,v}$ where  $a\in\E$ is monoidal and
Hermitian  and $u,v\in X$ is called the Tannaka groupoid of $\G$.
For fixed $u,v\in X$, we denote the set of all $a_{u,v}\in\T$ by
$\mathcal {T}_{u,v}(\G)$.
\end{defn}

\begin{thm} $\T$ is a compact groupoid.
\end{thm}
{\bf Proof} We define the product for the pairs of the form
$(a_{w,v},b_{u,w})\in\T^{(2)}$ by composition
$$(ab)_{u,v}^\pi=a_{w,v}^\pi\circ b_{u,w}^\pi.$$ This is clearly
an associative partial operation on $\T$.

It is easy to check that if $a,b\in\E$ are monoidal and Hermitian,
then so is $ab$. Indeed
\begin{align*}
(ab)_{u,v}^{\pi_1\oo\pi_2} &=a_{w,v}^{\pi_1\oo\pi_2}\circ
b_{u,w}^{\pi_1\oo\pi_2}=(a_{w,v}^{\pi_1}\oo
a_{w,v}^{\pi_2})\circ(b_{u,w}^{\pi_1}\oo
b_{u,w}^{\pi_2})\\
&=(a_{w,v}^{\pi_1}\circ b_{u,w}^{\pi_1})\oo (a_{w,v}^{\pi_2}\circ
b_{u,w}^{\pi_2})=(ab)_{u,v}^{\pi_1}\oo (ab)_{u,v}^{\pi_2}.
\end{align*}

For each $\pi\in\R$ let $\check\pi\in\R$ be its adjoint \repn ,
and put
$$(a_{u,v}^{-1})^\pi:={}^ta_{u,v}^{\check\pi}\quad (u,v\in
X,\pi\in\R).$$ For each $u\in X$ define $\eps_u:\mathcal
H_u^{\check\pi}\bigotimes\mathcal H_u^\pi\to\mathbb C$ by
$$\eps_u(\eta\oo\xi)=<\eta,\xi>\quad(\eta\in \mathcal H_u^{\check\pi}=(\mathcal
H_u^\pi)^{*},\xi\in\mathcal H_u^{\pi}),$$ then we claim that
$\eps\in Mor (\check\pi\oo\pi,tr)$. Indeed for each $x\in\G$ and
$\eta\in \mathcal H_{s(x)}^{\check\pi},\xi\in\mathcal
H_{s(x)}^{\pi}$ we have
\begin{align*}
\eps_{r(x)}\check\pi\oo\pi(x)(\eta\oo\xi)
&=\eps_{r(x)}(\check\pi(x)\eta\oo\pi(x)\xi)=<\check\pi(x)\eta,\pi(x)\xi>\\
&=<{}^t\pi(x)^{-1}\eta,\pi(x)\xi>=<\eta,\xi>\\
&=\eps_{s(x)}(\eta\oo\xi)=tr(x)\eps_{s(x)}(\eta\oo\xi).
\end{align*}
Therefore for each $u,v\in X$ and $a\in\E$ we have $\eps_v
a_{u,v}^{\check\pi\oo\pi}=a_{u,v}^{tr}\eps_u$. In particular for
each $a_{u,v}\in \T$, $\eta\in \mathcal H_{u}^{\check\pi}$, and
$\xi\in\mathcal H_{v}^{\pi}$ we have
\begin{align*}
<a_{u,v}^{\check\pi}(\eta),a_{u,v}^\pi(\xi)>&=\eps_v(
a_{u,v}^{\check\pi}(\eta)\oo
a_{u,v}^\pi(\xi))\\
&=\eps_v (a_{u,v}^{\check\pi\oo\pi}(\eta\oo
\xi))=a_{u,v}^{tr}\eps_u(\eta\oo\xi)=<\eta,\xi>.
\end{align*}

Put $\eta=b_{v,w}^{\check\pi}(\zeta)$ with $\zeta\in \mathcal
H_v^{\check\pi}$, then
$$<a_{u,v}^{\check\pi}(b_{v,w}^{\check\pi}(\zeta)),a_{u,v}^\pi(\xi)>
=<b_{v,w}^{\check\pi}(\zeta),\xi>,$$ for each $\zeta,\xi$ as
above. Hence, changing $\pi$ to $\check\pi$, we get
$${}^t a_{u,v}^{\check\pi}\circ a_{u,v}^{\pi}\circ b_{v,w}^{\pi}=b_{v,w}^{\pi},$$
that is $a_{v,u}^{-1}a_{u,v} b_{v,w}=b_{v,w}$. Similarly
$b_{w,v}a_{v,u}^{-1}a_{u,v} =b_{w,v}$. This shows that $\T$ is a
groupoid.

Next  we show that $\T$ is a closed subset of a compact groupoid.
Recall that isotropy groups $\G_u^u$ are compact groups and the
restriction of the invariant measure $d\lambda_u$ to $\G_u^u$ is a
left (and so right) Haar measure. For each $\pi\in\R$ and $u\in
X$, let $g_u:\Hu\bigotimes\Hu\to\mathbb C$ be defined by
$$ g_u(\xi,\eta)=\int_{G_u^u} <\pi(x)\xi,\eta>d\lambda_u^u(x)\quad(\xi,\eta\in\Hu).$$
Also  define $h_u:\overline{\Hu}\bigotimes\Hu\to\mathbb C$ by
$h_u(\xi,\eta)=g_u(\bar\xi,\eta)$. We claim that $h\in
Mor(\bar\pi\oo\pi,tr)$. Indeed for each $\xi,\eta\in\Hu$ and
$x\in\G$ we have
\begin{align*}
h_{r(x)}(\bar\pi\oo\pi)(x)(\xi\oo\eta)&=h_{r(x)}
(\bar\pi(x)\xi\oo\pi(x)\eta)=g_{r(x)}(\pi(x)\bar\xi\oo\pi(x)\eta)\\
&=\int<\pi(y)\pi(x)\bar\xi,\pi(x)\eta>d\lambda_{r(x)}^{r(x)}(y)\\
&=\int<\pi(x^{-1}yx)\bar\xi,\eta>d\lambda_{r(x)}^{r(x)}(y)\\
&=\int<\pi(y)\bar\xi,\eta>d\lambda_{s(x)}^{s(x)}(y)\\
&=g_{s(x)}(\bar\xi\oo\eta)=h_{s(x)}tr(x)(\xi\oo\eta).
\end{align*}
Therefore, for each $u,v\in X$ and $a\in \E$ we have $h_v
a_{u,v}^{\bar\pi\oo\pi}=a_{u,v}^{tr}h_u$. In particular for each
$a_{u,v}\in\T$ using monoidal property we get
$h_v(\overline{a_{u,v}^\pi(\xi)},
a_{u,v}^\pi(\eta))=h_u(\bar\xi,\eta)$, that is
$g_v(a_{u,v}^\pi(\xi), a_{u,v}^\pi(\eta))=g_u(\xi,\eta)$, for each
$\xi,\eta\in\Hu$. Now we can view $g_u$ and $g_v$ as new inner
products on $\Hu$ and $\Hv$, respectively, and look at the unitary
elements in $\mathbb B(\Hu,\Hv)$, then the above relation is just
to say that
$a_{u,v}^\pi\in\mathcal{U}(\mathcal{B}((\Hu,g_u),(\Hv,g_v)))$,
whence
$$\T\subseteq\prod_{\pi,u,v} \mathcal{U}(\mathcal{B}((\Hu,g_u),(\Hv,g_v))),$$
a product of compact groupoids.  The fact that $\T$ is a closed
subset of this groupoid follows immediately from the definition of
the topology on $\E$.\qed

\vspace{.3 cm} Now let's consider the natural transformations
$\mathcal T_x\in\E$, $x\in\G$. It is clear that for each $x\in\G$,
$\mathcal T_x\in\T$ and
$$\mathcal T_{xy}=\mathcal T_x\mathcal T_y\quad(x,y\in\G^{(2)}).$$
In particular the image of $\G$ under $\mathcal T$ is a
subgroupoid of $\T$. We identify $\G$ with its  image in $\T$. For
each $u,v\in X$, let $\mathcal
T_{u,v}:\G_u^v\to\mathcal{T}_{u,v}(\G)$ be defined by $\mathcal
T_{u,v}(x)(\pi)=\pi(x)\quad(x\in\G_u^v)$. Also we can define two
adjoint maps
$$\mathcal {T}^*: \mathcal {R}ep(\T)\to\R$$
by $$\mathcal {T}^*(\Pi)(x)=\Pi(\mathcal
T_x)\quad(x\in\G,\Pi\in\mathcal {R}ep(\T)),$$ and
$$\mathcal {T}^*: \mathcal E(\T)\to\mathcal E(\G)$$
by $$\mathcal {T}^*(f)(x)=f(\mathcal T_x)\quad(x\in\G,f\in\mathcal
E(\T)).$$

\begin{lem}
The restriction map
$$\mathcal {T}^*: \mathcal {R}ep(\T)\to\R$$
is a bundle isomorphism.
\end{lem}
{\bf Proof} We define the extension bundle map $\mathfrak
{E}:\R\to\mathcal {R}ep(\T)$  as follows. Given $u,v\in X$,
$a_{u,v}\in\T$, and $\pi\in\R$, the map $P_\pi: a_{u,v}\mapsto
a_{u,v}^\pi$ is a \repn of $\T$ on $\HH$ and we have the
commutative triangles

\begin{equation*}
\begin{CD}
\G_{u,v}@>{\pi}>>\mathcal B(\Hu,\Hv)\\
@V{\mathcal {T}_{u,v}}VV    @AA{P_{\pi}}A\\
\mathcal {T}_{u,v}(\G)@>>{=}>\mathcal {T}_{u,v}(\G)
\end{CD}
\end{equation*}

Therefore $\mathcal{T}^*(P_\pi)=\pi$. We  put $\mathfrak
E(\pi)=P_\pi$. If $h\in Mor_{\G}(\pi_1,\pi_2)$ then clearly $h\in
Mor_{\T}(P_{\pi_1},P_{\pi_2})$. Also it is easy to check that
$\mathfrak E$ preserves direct sums, tensor products, and
conjugation of \repns . Moreover the above commutative triangle
shows that if $\pi$ is irreducible, then so is $\mathfrak(\pi)$.
Hence $Im(\mathfrak E)$ is a closed subset of $\T$ in the sense of
Definition 3.7 of [A2]. It also separates the points of $\T$.
Indeed If $a_{u,v}$ and $b_{w,z}$ are distinct elements of $\T$,
there is a \repn $\pi\in\R$ such that $a_{u,v}^\pi\neq
b_{w,z}^\pi$, which means that $P_\pi$ separates $a_{u,v}$ and
$b_{w,z}$. By [A2, proposition 3.8], $\mathfrak E$ is surjective.
Now $\mathcal T^*\circ\mathfrak E=id$, so $\mathcal T^*$ is a
bundle isomorphism.\qed

\vspace{.3 cm} Now let $\mathcal E(\G)$ and $\mathcal E(\T)$ be
the  representation bundles of $\G$ and $\T$, respectively.

\begin{lem}
The restriction map
$$\mathcal {T}^*: \mathcal {R}ep(\T)\to\R$$
is a bundle isomorphism.
\end{lem}
{\bf Proof} We  define the extension bundle map $\mathfrak
{E}:\R\to\mathcal {R}ep(\T)$ as follows. Given $u,v\in X$, by
Proposition 2.2 of [A2], any $f\in\Euv$ has a unique \repn in the
form
$$f=\sum_{\pi\in\hat\G} d_u^\pi Tr\big(g(\pi))\pi(.)\big),$$
where  $g=\mathfrak F_{u,v}(f)\in \sum_{\pi\in\hat\G} \mathcal
B(\Hv,\Hu)$. Define $\mathfrak E_{u,v}(f)$ on $\mathcal
T_{u,v}(\G)$ by
$$\mathfrak E_{u,v}(f)(a_{u,v})=\sum_{\pi\in\hat\G} d_u^\pi Tr\big(g(\pi))
P_\pi(a_{u,v})\big)\quad (a_{u,v}\in \mathcal T_{u,v}(\G)).$$ By
above  lemma, $\mathfrak E$ is injective, so $\mathcal T^*$ is
bijective, as we have $\mathcal T^*\circ\mathfrak E=id$.\qed

\begin{lem}
For each $f\in C(\T)$ and $u,v\in X$,
$$\int_{\T_u^v} f(t)d\tilde\lambda_u^v(t)=\int_{\G_u^v} f(\pi(x))d\lambda_u^v(x).$$
\end{lem}
{\bf Proof} By Lemma 3.9 of [A1], it is  enough to prove this for
$f\in\mathcal E(\T)$. As in the proof of the above lemma we may
represent $f$ on $\T_u^v$ as
$$f(t)=\sum_{\pi\in\hat\G} d_u^\pi Tr\big(g(\pi))P_\pi(t)\big)\quad (t\in \T_u^v),$$
where $g=\mathfrak F(\mathcal T^*(f))$. In particular, for each
$x\in\G_u^v$,
$$f(\mathcal T_x))=\sum_{\pi\in\hat\G} d_u^\pi Tr\big(g(\pi))P_\pi(\mathcal T_x)\big)
=\sum_{\pi\in\hat\G} d_u^\pi Tr\big(g(\pi))\pi(x)\big).$$ By
Proposition 3.2 $(iii)$ of [A2], we have

\begin{equation*}
\int Tr\big(g(\pi)\pi(x)\big)d\lambda_u^v(x)=
\begin{cases}
\, g(tr) & \text{if} \,\, \pi=tr,\\
\,0& \text{otherwise},
\end{cases}
\end{equation*}
where $tr$ is the trivial \repn , and similarly

\begin{equation*}
\int Tr\big(g(\pi)P_\pi(t)\big)d\tilde\lambda_u^v(t)=
\begin{cases}
\, g(tr) & \text{if} \,\, \pi=tr,\\
\,0& \text{otherwise},
\end{cases}
\end{equation*}
hence the result.\qed

\vspace{.3 cm} Now we are ready to prove the main result of these
series of papers, the {\it Tannaka-Krein duality theorem} for
compact groupoids.

\begin{thm}{\bf (Tannaka-Krein Duality Theorem)}
Any compact groupoid is isomorphic to its Tannaka groupoid.
\end{thm}
{\bf Proof} Let $\G$ be  a compact groupoid. We show that
$\mathcal T:\G\to\T$ is an isomorphism of topological groupoids.
The injectivity of $\mathcal T$ follows from the Peter-Weyl
theorem [A1, theorem 3.13]. For the surjectivity, assume on the
contrary that $Im(\mathcal T)$ is a proper subset of $\T$. This is
a closed subset. Let $f\in C(\T)$ be a positive function such that
$supp(f)$ is contained in the complement of $Im(\mathcal T)$. Then
from the two integrals in above lemma, the one on the right hand
side is $0$, where as the one on the left hand side is strictly
positive, a contradiction. \qed


\end{document}